\documentclass[11pt]{article}
\usepackage{amsfonts,amssymb,latexsym,amscd, theorem,epsfig,graphics}

\newtheorem{theorem}{Theorem}

\begin{document}

\def\startproof{{\bf {\medskip}{\noindent}Proof: }}
\def\definition{{\bf {\bigskip}{\noindent}Definition: }}
\def\endproof{$\spadesuit$  \newline}

\def\bp{|\underline{\overline{\times}}|}
\def\A{\mbox{\boldmath{$A$}}}%
\def\B{\mbox{\boldmath{$B$}}}%
\def\C{\mbox{\boldmath{$C$}}}%
\def\D{\mbox{\boldmath{$D$}}}%
\def\E{\mbox{\boldmath{$E$}}}%
\def\F{\mbox{\boldmath{$F$}}}%
\def\G{\mbox{\boldmath{$G$}}}%
\def\H{\mbox{\boldmath{$H$}}}%
\def\I{\mbox{\boldmath{$I$}}}%
\def\J{\mbox{\boldmath{$J$}}}%
\def\K{\mbox{\boldmath{$K$}}}%
\def\L{\mbox{\boldmath{$L$}}}%
\def\M{\mbox{\boldmath{$M$}}}%
\def\N{\mbox{\boldmath{$N$}}}%
\def\O{\mbox{\boldmath{$O$}}}%
\def\P{\mbox{\boldmath{$P$}}}%
\def\Q{\mbox{\boldmath{$Q$}}}%
\def\R{\mbox{\boldmath{$R$}}}%
\def\T{\mbox{\boldmath{$T$}}}%
\def\U{\mbox{\boldmath{$U$}}}%
\def\V{\mbox{\boldmath{$V$}}}%
\def\W{\mbox{\boldmath{$W$}}}%
\def\X{\mbox{\boldmath{$X$}}}%
\def\Y{\mbox{\boldmath{$Y$}}}%
\def\Z{\mbox{\boldmath{$Z$}}}%

\title {Elementary Surprises in Projective Geometry}
\author{Richard Evan Schwartz\thanks{Supported by 
N.S.F. Research Grant DMS-0072607.}\hskip 3 mm  and
Serge Tabachnikov\thanks{Supported by 
N.S.F. Research Grant DMS-0555803. Many thanks to MPIM-Bonn for its hospitality.}}
\date{}
\maketitle

\begin{figure}[hbtp]
\centering
\includegraphics[width=3in]{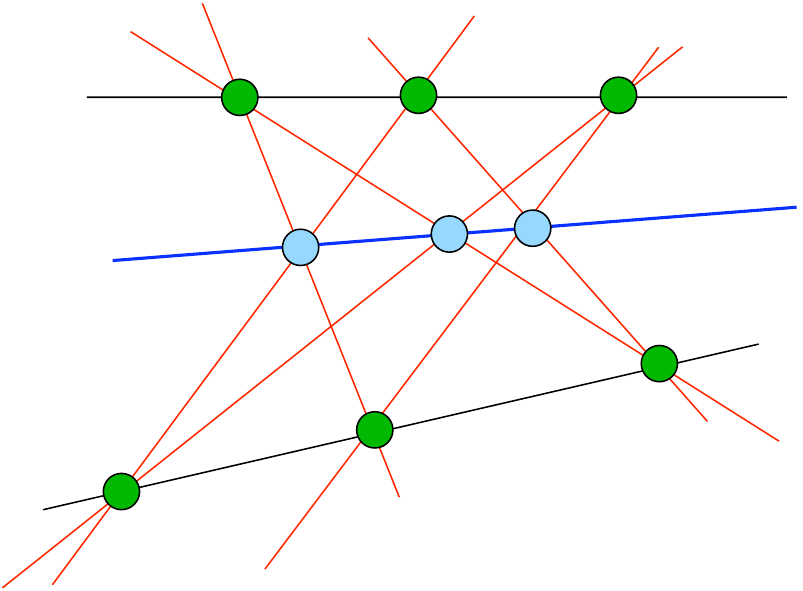}
\caption{Pappus's Theorem}
\label{pap}
\end{figure}

The classical theorems in projective geometry involve constructions
based on points and straight lines.  A general feature of these theorems
is that a surprising coincidence awaits the reader who makes the
construction.  One example of this is Pappus's theorem.
One starts with $6$ points, $3$ of which are contained on
one line and $3$ of which are contained on another.
Drawing the additional lines shown in Figure \ref{pap}, one sees that
the $3$ middle (blue) points are also contained on a line.

\begin{figure}[hbtp]
\centering
\includegraphics[width=5in]{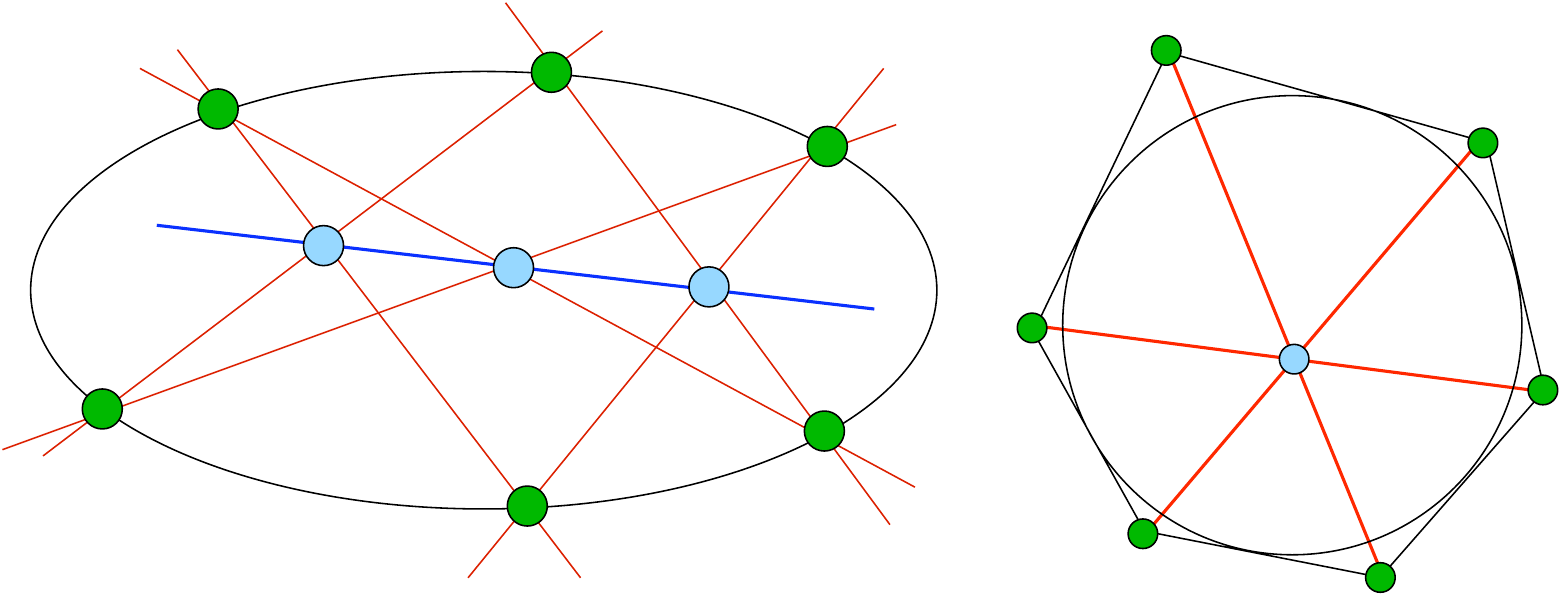}
\caption{Pascal's Theorem and Brian\c{c}on's Theorem}
\label{pas}
\end{figure}

Pappus's Theorem goes back about $1700$ years.  In $1639$, Blaise Pascal
discovered a generalization of Pappus's Theorem.  In Pascal's
Theorem, the $6$ green points are contained in a conic section, as
shown on the left hand side of Figure \ref{pas}.

One recovers Pappus's Theorem as a kind of limit, as the conic section
stretches out and degenerates into a pair of straight lines.

Another closely related theorem is Brian\c{c}on's Theorem.  This
time, the $6$ green points are the vertices of a hexagon that
is circumscribed about a conic section, as shown on the right hand side of
Figure \ref{pas}, and the surprise is that
the $3$ thickly drawn diagonals intersect in a point.
 Though
Brian\c{c}on discovered this result about $200$ years after
Pascal's theorem, the two results are in fact equivalent for the
well-known reason we will discuss below.

The purpose of this article is to discuss some apparently new
theorems in projective geometry that are similar in spirit
to Pascal's Theorem and Brian\c{c}on's Theorem.  
One can think of all the results we discuss as statements
about lines and points in the ordinary Euclidean plane,
but setting the theorems in the {\it projective plane\/}
enhances them.
\newline
\newline
{\bf The Basics of Projective Geometry:\/}
Recall that the projective plane $\P$ is defined as the space of
lines through the origin in $\R^3$.  A point in $\P$ can be described
by {\it homogeneous coordinates\/}
$(x:y:z)$, not all zero, corresponding to the line containing
the vector $(x,y,z)$.  Of course, the two
triples $(x:y:z)$ and $(ax:ay:az)$ describe
the same point in $\P$ as long as $a \not = 0$. One says that $\P$ is the {\it projectivization\/} of $\R^3$.

  A {\it line\/} in
the projective plane is defined as a set of lines through the origin in $\R^3$ that lie in a plane through the origin.
Any linear isomorphism of $\R^3$ -- i.e., multiplication
by an invertible $3 \times 3$ matrix -- permutes the
lines and planes through the origin. In this way, a linear
isomorphism induces a mapping of $\P$ that carries
lines to lines.   These maps are called {\it projective
transformations\/}.  

One way to define a (non-degenerate) conic section in $\P$ is to say
that
\begin{itemize}
\item The set of points in $\P$ of the form $(x:y:z)$ such
that $z^2=x^2+y^2 \not = 0$ is a conic section.
\item Any other conic section is the image of the one
we just described under a projective transformation.
\end{itemize}

One frequently identifies $\R^2$ as the subset
of $\P$ corresponding to points $(x:y:1)$.
We will simply write $\R^2 \subset \P$.
The ordinary lines in $\R^2$ are subsets of
lines in $\P$.   The conic sections
intersect $\R^2$ in either ellipses,
hyperbolas, or parabolas.  One of the
beautiful things about projective
geometry  is that these three
kinds of curves are {\it the same\/} from
the point of view of the projective plane
and its symmetries.

The {\it dual plane\/} $\P^*$ is defined  to be
the set of planes through the origin in $\R^3$.
Every such plane is  the kernel of a linear function on $\R^3$, and this linear function is determined by the plane up to a non-zero factor. Hence $\P^*$ is the projectivization of the dual space $(\R^3)^*$. If one wishes, one can identify $\R^3$ with $(\R^3)^*$ using the scalar product. 
One can also think of $\P^*$ as the space of lines in $\P$. 

\begin{figure}[hbtp]
\centering
\includegraphics[width=2.2in]{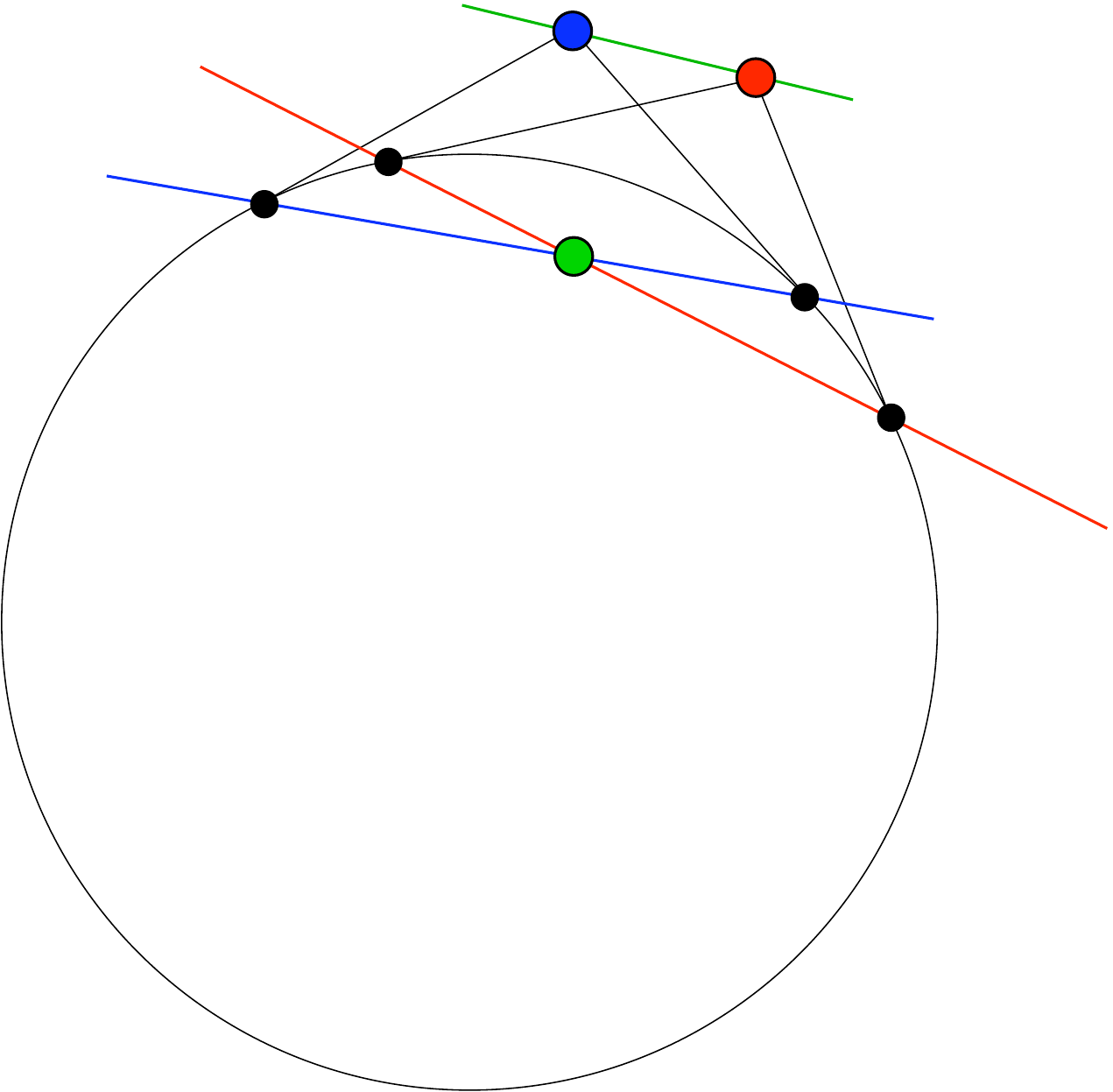}
\caption{Projective duality}
\label{dual}
\end{figure}

Given a point $v$ in $\P$, the set $v^{\perp\/}$ of linear functions on $\R^3$, that vanish at $v$, determine a line in $\P^*$. The correspondence $v\mapsto v^{\perp}$ carries collinear points to concurrent lines; it is called the 
{\it projective duality\/}.  A projective duality takes 
points of $\P$ to lines of $\P^*$, and lines of $\P$ to points of $\P^*$. 
Of course, the same construction works in the opposite direction, from $\P^*$ to $\P$. 
Projective duality is an involution: applied twice, it yields the identity map.   Figure \ref{dual} illustrates
an example of a projective duality based on the unit circle:   the red line maps to the red point, the blue line maps to the blue point, and the green point maps to the green line.

Projective duality extends to smooth curves: the 1-parameter family of the tangent lines to a curve $\gamma$ in $\P$ is a 1-parameter family of points in $\P^*$, the dual curve $\gamma^*$. The curve dual to a conic section is again a conic section. Thus projective duality carries the vertices
of a polygon inscribed in a conic to the lines extending
the edges of a polygon circumscribed about a conic.

Projective duality  
takes an instance of Pascal's Theorem to an
instance of Brian\c{c}on's Theorem, and
vice versa.    This becomes clear if
one looks at the objects involved. The input
of Pascal's theorem is an inscribed hexagon and
the output is $3$ collinear points.
The input of Brian\c{c}on's theorem is a
superscribed hexagon and the output is 
$3$ coincident lines.
\newline
\newline
{\bf Polygons:\/} 
Like Pascal's Theorem and Brian\c{c}on's Theorem, our results all
involve polygons. A polygon $P$ in $\P$ is a cyclically ordered collection
$\{p_1,...,p_n\}$ of points, its vertices. A polygon has sides: the cyclically ordered collection
$\{l_1,...,l_n\}$ of lines in $\P$ where $l_i=\overline{p_ip_{i+1}}$ for all $i$. Of course, the indices are taken mod $n$. The {\it dual polygon} $P^*$ is the polygon in $\P^*$ whose vertices are $\{l_1,...,l_n\}$; the sides of the dual polygon are $\{p_1,...,p_n\}$ (considered as lines in $\P^*$). The polygon dual to the dual is the original one: $(P^*)^*=P$.

Let ${\cal X}_n$ and ${\cal X}_n^*$ denote the sets of $n$-gons in $\P$ and $\P^*$, respectively. 
There is a natural
map $ T_k: {\cal X}_n  \to  {\cal X}_n^*.$ 
Given an $n$-gon $P=\{p_1,...,p_n\}$, we define
$T_k(P)$  as
$$
\{\overline{p_1p_{k+1}},\overline{p_2p_{k+2}},\ldots \overline{p_n,p_{k+n}}\}.
$$
That is, the vertices of $T_k(P)$ are the consecutive $k$-diagonals of $P$.
The map $T_k$ is an {\it involution\/}, meaning that $T_k^2$ is
the identity map.     
When $k=1$, the map $T_1$ carries a polygon to the dual one.

\begin{figure}[hbtp]
\centering
\includegraphics[width=2in]{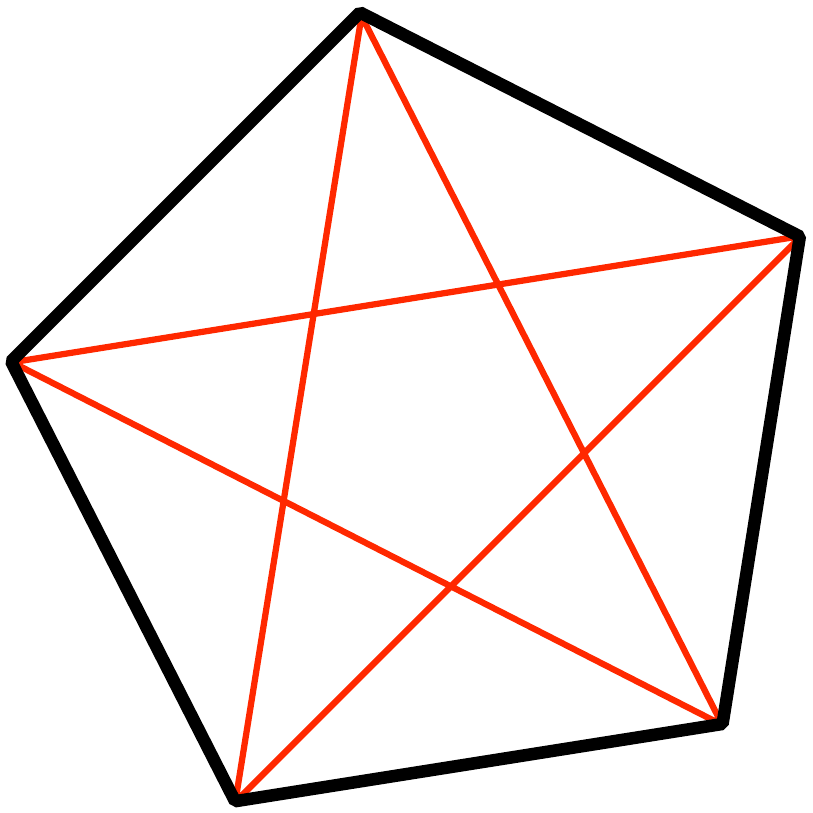}
\caption{The pentagram}
\label{pentag}
\end{figure}

Even when $a \not = b$, the map $T_{ab}=T_a \circ T_b$
carries ${\cal X}_n$ to ${\cal X}_n$ and ${\cal X}_n^*$ to ${\cal X}_n^*$.  We
have studied the dynamics of the {\it pentagram map} $T_{12}$
in detail in \cite{Sch1,Sch2,Sch3,OST1,OST2}, and the configuration theorems we present here are a byproduct of that study. (The map is so-called because of the resemblence, in the special case of
pentagons, to the famous mystical symbol having the same name. See Figure \ref{pentag}.)  We extend the notation: $T_{abc}=T_a \circ T_b\circ T_c$, and so on.

Now we are ready to present our configuration theorems.
\newline
\newline
{\bf The Theorems:\/}
To save words,
we say that an {\it inscribed polygon\/} is a polygon whose
vertices are contained in a conic section.   Likewise, we say that a {\it circumscribed polygon\/} is a polygon whose sides are tangent to a conic. Projective duality carries inscribed polygons to circumscribed ones and vice versa.  
We say that two polygons, $P$ in $\P$ and $Q$ in $\P^*$, are {\it equivalent\/} 
if there is a projective transformation $\P\to\P^*$ that takes $P$ to $Q$.  In this case, we write $P \sim Q$.  
By {\it projective transformation $\P \to \P^*$\/} we mean a map that
 is induced by a linear map $\R^3\to (\R^3)^*$. 

\begin{figure}[hbtp]
\centering
\includegraphics[width=3.7in]{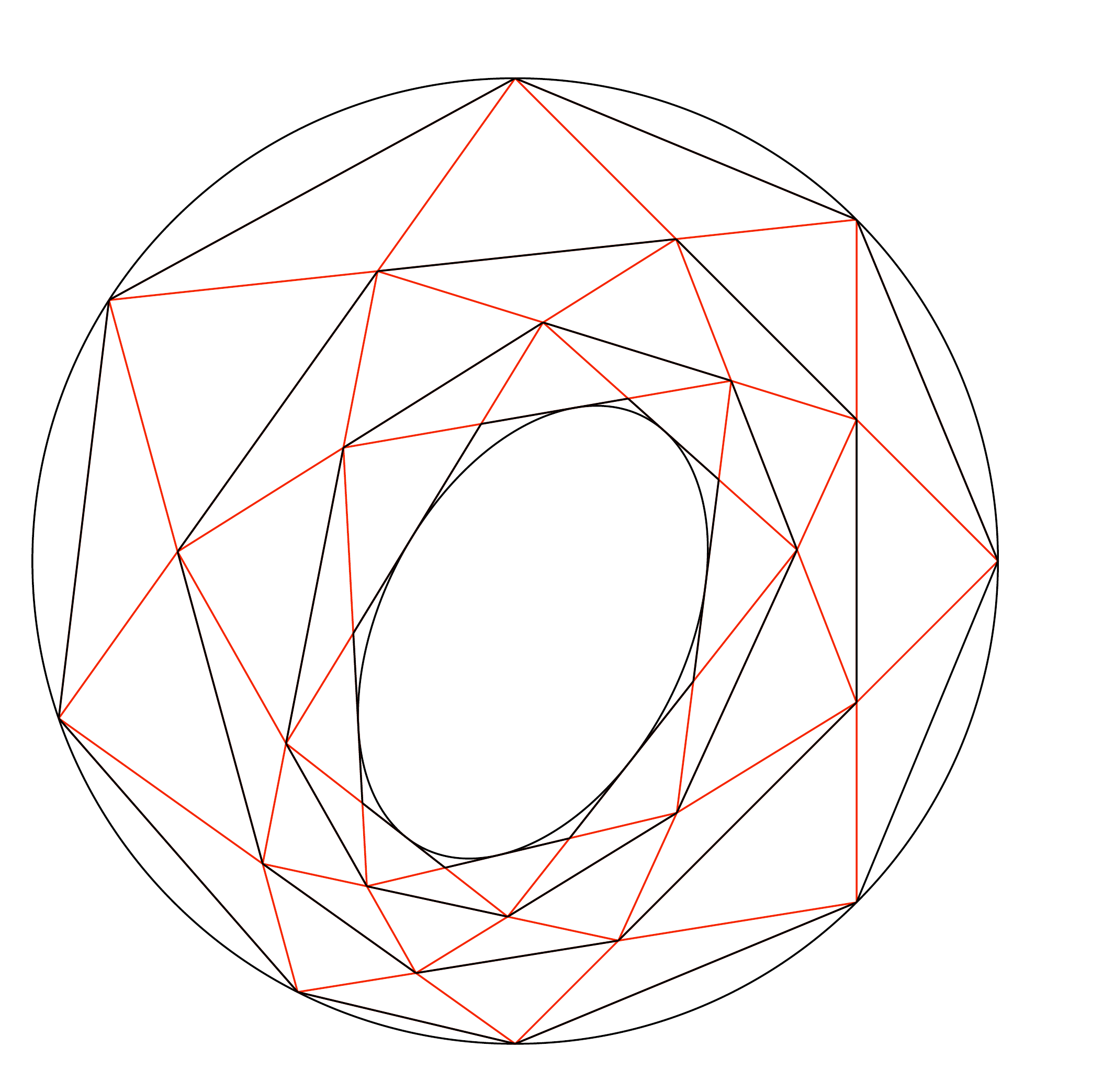}
\caption{If $P$ is an inscribed octagon then $P \sim T_{21212}(P)$}
\label{thm1}
\end{figure}

\begin{theorem} \label{main}
The following is true.
\begin{itemize}
\item If $P$ is an inscribed $6$-gon, then $P \sim T_2(P)$.
\item If $P$ is an inscribed $7$-gon, then $P \sim T_{212}(P)$.
\item If $P$ is an inscribed $8$-gon, then $P \sim T_{21212}(P)$.
\end{itemize}
\end{theorem}
Figure \ref{thm1} illustrates\footnote{Our Java applet does a much
better job illustrating these results.  To play with it online, see
http://www.math.brown.edu/$\sim$res/Java/Special/Main.html.}
 the third of these results.
The outer octagon $P$ is inscribed in a conic and the innermost
octagon $T_{121212}(P)=(T_{21212}(P))^*$ is circumscribed about a conic.

The reader might wonder if our three results are the beginning of an
infinite pattern.  Alas, it is not true that
$P$ and $T_{2121212}(P)$ are equivalent when $P$ is in inscribed $9$-gon,
and the predicted result fails for larger $n$ as well.  However,
we do have a similar result for $n=9,12$.

\begin{theorem} \label{9gon}
If $P$ is a circumscribed $9$-gon, then $P \sim T_{313}(P)$.
\end{theorem}

\begin{theorem} \label{12gon}
If $P$ is an inscribed $12$-gon, then
$P \sim T_{3434343}(P)$.
\end{theorem}

Even though all conics are projectively equivalent, it is not true
that all $n$-gons are projectively equivalent.   For instance,
the space of inscribed $n$-gons, modulo projective
equivalence, is $n-3$ dimensional.
We mention this because our last collection of results all
make weaker statements to the effect that the ``final
polygon''' is cicumscribed but not necessarily equivalent or projectively dual to
the ``initial polygon''.

\begin{theorem}
\label{weak}
The following is true.
\begin{itemize}
\item If $P$ is an inscribed $8$-gon, then $T_{3}(P)$ is circumscribed.
\item If $P$ is an inscribed $10$-gon, then $T_{313}(P)$ is circumscribed.
\item (*) If $P$ is an inscribed $12$-gon, then $T_{31313}(P)$ is circumscribed.
\end{itemize}
\end{theorem}

\begin{figure}[hbtp]
\centering
\includegraphics[width=3.7in]{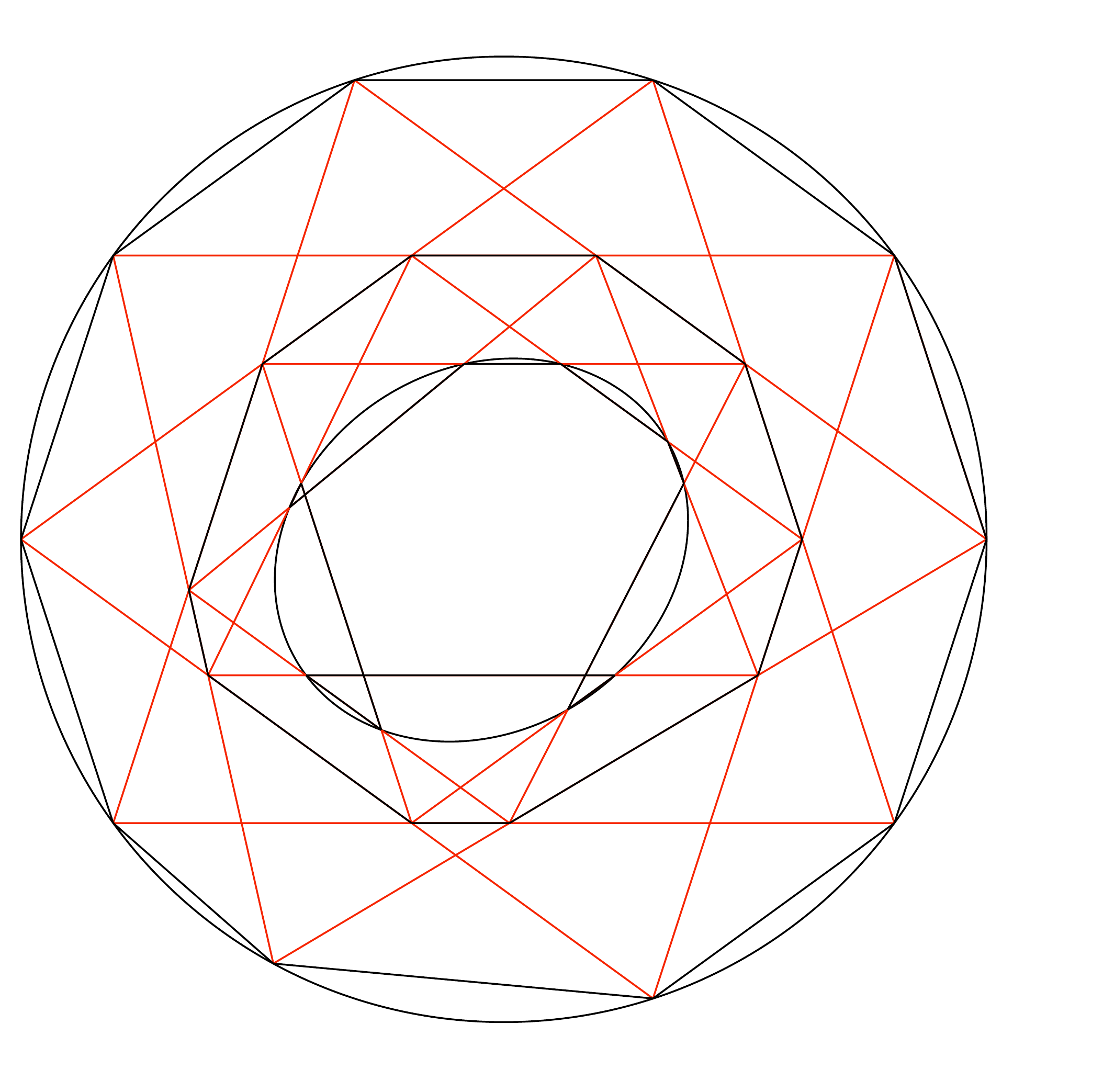}
\caption{ If $P$ is an inscribed decagon  then $T_{1313}(P)$ is also inscribed}
\label{thm2}
\end{figure}

We have starred the third result because we don't yet have a proof
for this one.  Figure \ref{thm2} illustrates the second of these results.
The formulation in Figure \ref{thm2} is easily seen to be equivalent
to the formulation given in Theorem \ref{weak}.
Looking carefully, we see that $T_{1313}(P)$ is not even convex.
(Even though the map $T_{13}$ is well defined on the
subset of convex polygons, it is not true that $T_{13}$
preserves this set.)
So, even though $T_{1313}(P)$ is inscribed, it is not projectively
equivalent to $P$ nor to its dual $P^*$. 
One might wonder if this result is part of an infinite pattern,
but once again the pattern stops after $n=12$.
\newline
\newline
{\bf Discovery and Proof:\/}
We discovered these results through computer experimentation.   We have
been studying the dynamics of the pentagram map $T_{12}$ on general polygons,
and we asked ourselves whether we could expect any special
relations when the intial polygon was either inscribed or circumscribed.

We initially found the $7$-gon result mentioned above. Then V. Zakharevich, 
a participant of the Penn State REU (Research Experience for Undergraduates) program in 2009, found Theorem \ref{9gon}.
Encouraged by this good luck, we made a more extensive computer search that
turned up the remaining results.  We think that the list above is
exhaustive, in the sense that there aren't any other surprises to be
found by applying some combination of diagonal maps to inscribed
or superscribed polygons.  In particular, we don't think that surprises
like the ones we found exist for $N$-gons with $N>12$.

The reader might wonder how we prove the results above.   In several
of the cases, we found some nice geometric proofs which we will
describe in a longer version of this article.  With one exception,
we found uninspiring algebraic proofs for the remaining cases.
Here is a brief description of these algebraic proofs.
First, we use symmetries of the projective plane
to reduce to the case when  
the vertices of $P$ lie on the parabola $y=x^2$.  We represent
vertices of $P$ in homogeneous coordinates in the form $(t:t^2:1)$.
Computing the maps $T_{k}(P)$ involves taking some cross products of
the vectors $(t,t^2,1)$ in $\R^3$.  At the end of the construction, our claims
about the final polygon boil down to equalities between
determinants of various $3 \times 3$ matrices made from
the vectors we generate.  We then check these identities
symbolically.

This approach has served to prove all but one of our
results:  the starred case of Theorem \ref{weak}.  The intensive symbolic
manipulation required for this case is currently beyond
what we can manage in Mathematica.   We don't know for
sure -- because we can't actually make the computation--
 but we think that the relevant polynomials (in $9$
variables) would have more than a trillion terms.
Naturally, we hope for
some clever cancellations that we haven't yet been
able to find.

We hope to find nice proofs for all the results above, but
so far this has eluded us.  Perhaps the interested reader
will be inspired to look for nice proofs.   We also hope that
these results point out some of the beauty of the dynamical
systems defined by these iterated diagonal maps.  Finally,
we wonder if the isolated results we have found are
part of an infinite pattern.  We don't have an opinion
one way or the other whether this is the case, but we
think that something interesting must be going on.
\newline
\newline
{\bf Additional Remarks:\/} 
In this concluding section, we relate our results to some other
classical constructions in projective geometry, and also give some
additional perspective on them.
\newline
\newline
1).  Let us say a few words about pentagons. The following is true:
\begin{itemize}
\item Every pentagon is inscribed in a conic and circumscribed about a conic. 
\item Every pentagon is projectively equivalent to its dual. 
\item The pentagram map  is the identity for every pentagon: $T_{12}(P)=P$. 
\end{itemize}
We do not want to deprive the reader from the pleasure of discovering proofs to the latter two claims (in case of difficulty, see \cite{FT} and \cite{Sch1}). Therefore one may add the following to Theorem \ref{main}: {\it If $P$ is a  $5$-gon, then $P\sim T_2(P)$}.

Related to the second item above, is the notion of a {\it self-polar} spherical polygon. Let $p_1,\dots,p_5$ be the vertices of a spherical pentagon. The pentagon is called self-polar if, for all $i=1,\dots,5$,  choosing $p_i$ as a pole, the points $p_{i+2}$ and $p_{i+3}$ both lie on the equator.  C. F. Gauss studied the geometry of such pentagons in a posthumously published work  {\it Pentagramma Mirificum}.

2). The formulations of Theorems \ref{main}-\ref{weak} are similar: if $P$ is inscribed, or circumscribed,  then $T_w(P)$ is projectively equivalent to $P$ (or is circumscribed). Here $w$ is a word in symbols $1,2,3,4$, that varies from statement to statement, but in each case, $w$ is {\it palindromic}: it is the same whether we read it left to right or right to left. This implies that, in each case, the transformation $T_w$ is an involution: $T_w \circ T_w=Id$.
\smallskip

3). The statement of Theorem \ref{main} can be rephrased as follows: {\it if $P$ is an inscribed heptagon then $T_2(P)$ and $T_{12}(P)$ are projectively equivalent}. That is, the heptagon $Q=T_2(P)$ is equivalent to its projective dual $Q^*$. In fact, every projectively self-dual heptagon is obtained this way. 

Similarly, Theorem \ref{9gon} states: {\it if $P$ is a circumscribed nonagon then $T_3(P)$ and $T_{13}(P)$ are projectively equivalent}, and hence $Q=T_3(P)$ is projectively self-dual. Once again, every projectively self-dual nonagon is obtained this way.

For odd $n$, the space of projectively self-dual $n$-gons in the projective plane, considered up to projective equivalence, is $n-3$-dimensional, see \cite{FT} (compare with $2n-8$, the dimension of projective equivalence classes of all $n$-gons). The space of inscribed (or circumscribed) $n$-gons, considered up to projective equivalence of the conic, is also $n-3$-dimensional. Thus, for $n=7$ and $n=9$, we have explicit bijections between these spaces.
\smallskip

4). One may cyclically relabel the vertices of a polygon to deduce apparently new configuration theorems from Theorems \ref{main}-\ref{weak}. Let us illustrate this by an example. Rephrase the last statement of Theorem \ref{weak} as follows: {\it If $P$ is an inscribed dodecagon then $T_{131313}(P)$ is also inscribed}. Now relabel the vertices as follows: $\sigma(i)=5i$ mod 12 (note that $\sigma$ is an involution). The map $T_3$ is conjugated by $\sigma$ as follows: 
$$
i\mapsto 5i\mapsto 5i+3\mapsto 5(5i+3)=i+3\ \ {\rm mod}\ 12, 
$$
that is, the map is $T_3$ again, and the map $T_1$ becomes
$$
i\mapsto 5i\mapsto 5i+1\mapsto 5(5i+1)=i+5\ \ {\rm mod}\ 12,
$$
that is, the map is $T_5$. We arrive at the statement: {\it If $P$ is an inscribed dodecagon then $T_{535353}(P)$ is also inscribed}.
Our java applet, cited above, shows pictures of this.
\smallskip

5). Theorem \ref{weak} appears to be a relative of a theorem in \cite{Sch3}: {\it Let $P$ be a $4n$-gon whose odd sides pass through one fixed point and whose even sides pass through another fixed point. Then the $(2n-2)$nd iterate of the pentagram map $T_{12}$ transforms $P$ to a polygon whose odd vertices lie on one fixed line and whose even vertices lie on another fixed line.} Note that a pair of lines is a degenerate conic section. Note also that the dual polygon, $Q=T_1(P)$ is also inscribed into a pair of lines. Thus we have an equivalent formulation: {\it If $Q$ is a $4n$-gon inscribed into a degenerate conic then $(T_1T_2)^{2n-2}T_1(Q)$ is also inscribed into a degenerate conic}.

We wonder if this result is a degenerate case of a more general theorem, much in the same way that
Pappus's theorem is a degenerate case of Pascal's theorem.

\end{document}